\documentclass[11pt,reqno]{amsart}
\usepackage{amsmath}
\usepackage{amssymb}
\usepackage{latexsym}
\usepackage{mathrsfs}

\def\squarebox#1{\hbox to #1{\hfill\vbox to #1{\vfill}}}

\theoremstyle{plain}

\newtheorem{thm}{Theorem}

\newtheorem{lem}{Lemma}

\newtheorem{prop}{Proposition}

\newtheorem{defn}{Definition}

\def\ii{{\bf i}}
\numberwithin{equation}{section}

\begin{document}
\def\C{{\mathbb C}}
\def\R{{\mathbb R}}
\def\N{{\mathbb N}}
\def\Z{{\mathbb Z}}
\def\T{{\mathbb T}}
\def\Q{{\mathbb Q}}
\def\SP{{\mathbb S}}
\def\d{{\partial}}
\def\mc{{\mathcal H}}
\def\1b{{\mathbb I}}
\def\tr{{\rm tr}\:}
\def\pv{\partial_x V}
\def\re{{\rm Re}\:}
\def\im{{\rm Im}\:}
%%%%%%%%%%%%%%%%%%%%%%%%%%%%%%%%%%%%%%%%%%%%%%%%%%%%%%%%%%%%%%%%%%%%%%
%%%%%%%%%%%%%%%%%%%%%%%%  Local definitions  %%%%%%%%%%%%%%%%%%%%%%%%%
%%%%%%%%%%%%%%%%%%%%%%%%%%%%%%%%%%%%%%%%%%%%%%%%%%%%%%%%%%%%%%%%%%%%%%
%

%\headheight = 25pt

%
%%%%%%%%%%%%%%%%%%%%%%%%%%%%%%%%%%%%%%%%%%%%%%%%%%%%%%%%%%%%%%%%%%%%%%
%%%%%%%%%%%%%%%%%  Uncomment to avoid overfull messages  %%%%%%%%%%%%%
%%%%%%%%%%%%%%%%%%%%%%%%%%%%%%%%%%%%%%%%%%%%%%%%%%%%%%%%%%%%%%%%%%%%%%
%
 \vbadness=10000
 \hbadness=10000
%
%%%%%%%%%%%%%%%%%%%%%%%%%%%%%%%%%%%%%%%%%%%%%%%%%%%%%%%%%%%%%%%%%%%%%%
%%%%%%%%%%%%%%%%%%%    General layout macros    %%%%%%%%%%%%%%%%%%%%%%
%%%%%%%%%%%%%%%%%%%%%%%%%%%%%%%%%%%%%%%%%%%%%%%%%%%%%%%%%%%%%%%%%%%%%%
%

%

% \renewcommand{\thetheorem}{\thesection.\arabic{theorem}}
% \renewcommand{\theproposition}{\thesection.\arabic{proposition}}
% \renewcommand{\thelemma}{\thesection.\arabic{lemma}}
% \renewcommand{\thedefinition}{\thesection.\arabic{definition}}
% \renewcommand{\thecorollary}{\thesection.\arabic{corollary}}
% \renewcommand{\theequation}{\thesection.\arabic{equation}}
% \renewcommand{\theremark}{\thesection.\arabic{remark}}

%
%%%%%%%%%%%%%%%%%%%%%%%%%%%%%%%%%%%%%%%%%%%%%%%%%%%%%%%%%%%%%%%%%%%%%%
%
    \setlength{\textwidth}{5.4 in}
%    \setlength{\textheight}{8 in}
%
%%%%%%%%%%%%%%%%%%%%%%%%%%%%%%%%%%%%%%%%%%%%%%%%%%%%%%%%%%%%%%%%%%%%%%
%%%%%%%%%%%%%%%%%%%%%%%%%%%%%%%%%%%%%%%%%%%%%%%%%%%%%%%%%%%%%%%%%%%%%%
%%%%%%%%%%%%%%%%%%%%%%%%%%%%%%%%%%%%%%%%%%%%%%%%%%%%%%%%%%%%%%%%%%%%%%
%%%%%%%%%%%%%%%%%%%%%%   End of preamble  %%%%%%%%%%%%%%%%%%%%%%%%%%%%
%%%%%%%%%%%%%%%%%%%%%%%%%%%%%%%%%%%%%%%%%%%%%%%%%%%%%%%%%%%%%%%%%%%%%%
%%%%%%%%%%%%%%%%%%%%%%%%%%%%%%%%%%%%%%%%%%%%%%%%%%%%%%%%%%%%%%%%%%%%%%
%%%%%%%%%%%%%%%%%%%%%%%%%%%%%%%%%%%%%%%%%%%%%%%%%%%%%%%%%%%%%%%%%%%%%%

%
% \baselineskip 18pt
%%%%%%%%%%%%%%%%%%%%%%%%%%%%%%%%%%%%%%%%%%%%%%%%%%%%%%%%%%%%%%%%%%%%%%
%%%%%%%%%%%%%%%% Blackboard boldface common symbols  %%%%%%%%%%%%%%%%%
%%%%%%%%%%%%%%%%%%%%%%%%%%%%%%%%%%%%%%%%%%%%%%%%%%%%%%%%%%%%%%%%%%%%%%
\def\R {{\mathbb{R}}}
\def\N {{\mathbb{N}}}
\def\C {{\mathbb{C}}}
\def\Z {{\mathbb{Z}}}
%%%%%%%%%%%%%%%%%%%%%%%%%%%%%%%%%%%%%%%%%%%%%%%%%%%%%%%%%%%%%%%%%%%%%%
%%%%%%%%%%%%%%%%%% varsymbols ... for readability %%%%%%%%%%%%%%%%%%%%
%%%%%%%%%%%%%%%%%%%%%%%%%%%%%%%%%%%%%%%%%%%%%%%%%%%%%%%%%%%%%%%%%%%%%%
\def\phi{\varphi}
\def\epsilon{\varepsilon}
\def\kappa{\varkappa}
\def\el{e^{-\lambda t}}
\def\bv {b_1(t,\xi)u}
\def\ia{t^2 a_3(t, \xi) \im (u \bar{u}'')}
\def\buu{\bar{u}''}
\def\bu{\bar{u}'}
\def\eT{e^{-\lambda T}}
\def\ii{{\bf i}}
\def\hh{\hat{x}}
\def\hhx{\hat{\xi}}
\def\mathscr{{\mathcal }}
%%%%%%%%%%%%%%%%%%%%%%%%%%%%%%%%%%%%%%%%%%%%%%%%%%%%%%%%%%%%%%%%%%%%%%
%%%%%%%%%%%%%%%%%%%%%%%%%%%%%%%%%%%%%%%%%%%%%%%%%%%%%%%%%%%%%%%%%%%%%%
%
%%%%%%%%%%%%%%%%%%%%%%%%%%%%%%%%%%%%%%%%%%%%%%%%%%%%%%%%%%%%%%%%%%%%%%
%%%%%%%%%%%%%%%%%%%%  math local macros %%%%%%%%%%%%%%%%%%%%%%%%%%%%%%
%%%%%%%%%%%%%%%%%%%%%%%%%%%%%%%%%%%%%%%%%%%%%%%%%%%%%%%%%%%%%%%%%%%%%%
%\def\tb#1{\|\kern -1.2pt | #1 \|\kern -1.2pt |} 
\def\Qed{\qed\par\medskip\noindent}
%%%%%%%%%%%%%%%%%%%%%%%%%%%%%%%%%%%%%%%%%%%%%%%%%%%%%%%%%%%%%%%%%%%%%%
%%%%%%%%%%%%%%%%%%%%%%%%%%%%%%%%%%%%%%%%%%%%%%%%%%%%%%%%%%%%%%%%%%%%%%
%
\title[Cauchy problem for hyperbolic operators ]{Cauchy problem for hyperbolic operators with triple characteristics of variable multiplicity}
\date{}
\author[E. Bernardi, A. Bove] {Enrico Bernardi, Antonio Bove}
\author[V. Petkov]{Vesselin Petkov}

\address{Dipartimento di Matematica per le Scienze Economiche e Sociali, Universit\`a di Bologna,
Viale Filopanti 5, 40126 Bologna, Italia}

\email{bernardi@economia.unibo.it}
\address{Dipartamento di Matematica,  Universit\`a di Bologna,\:
Piazza di Porta S. Donato 5, 
40126 Bologna, Italia}
\email{bove@dm.unibo.it}

\address{Universit\'e Bordeaux I, Institut de Math\'ematiques de Bordeaux, 351, 
Cours de la Lib\'eration, 33405  Talence, France}
\email{petkov@math.u-bordeaux1.fr}
\maketitle
\begin{abstract}
We study a class of third order hyperbolic operators $P$ in $G = \Omega \cap \{0 \leq t \leq T\},\: \Omega \subset \R^{n+1}$ with triple characteristics on $t = 0$. We consider the case when the fundamental matrix of the principal symbol for $t = 0$ has a couple of non vanishing real eigenvalues and $P$ is strictly hyperbolic for $t > 0.$ We prove that $P$ is strongly hyperbolic, that is the Cauchy problem for $P + Q$ is well posed in $G$ for any lower order terms $Q$.
\end{abstract}

\section{Introduction}
\renewcommand{\theequation}{\arabic{section}.\arabic{equation}}
\setcounter{equation}{0}
Consider a differential operator
 $$P(t, x, D_t, D_x) = \sum _{\alpha + |\beta| \leq m} c_{\alpha, \beta} (t, x) D_t^{\alpha} D_x^{\beta},\: D_t = -\ii\partial_t, D_{x_j} = -\ii\partial_{x_j}$$ 
 of order $m$ with $C^{\infty}$ coefficients $c_{\alpha, \beta}(t,x),\: t \in \R,\:x \in \R^n.$ Denote by 
$$p_m(t, x, \tau, \xi) = \sum_{\alpha + |\beta| = m} c_{\alpha, \beta} (t, x) \tau^{\alpha} \xi^{\beta}$$
 the principal symbol of $P$.
Let $\Omega \subset \R^{n + 1}$ be an open set and let 
$$\Omega_{\eta}^{-} = \Omega \cap \{t \leq \eta\},\Omega_{\eta}^{+} = \Omega \cap \{t \geq \eta\}, \:G =  \Omega \cap \{0 \leq t \leq T\}.$$
 We say that $P$ is hyperbolic with respect to $N_0 =(1,0,...,0)$ at $(t_0, x_0)$ if
$$ (i)\:\:\:\: p_m(t_0, x_0, N_0) \neq 0,$$ 
(ii)\:\:the equation 
\begin{equation} \label{eq:1.1}
p_m(t_0, x_0, \tau, \xi) = 0
\end{equation}
 with respect to $\tau $ has only real roots  $\tau = \lambda_j(t_0, x_0, \xi)$ for all $\xi \in \R^n.$
Set $P_m(t, x, D_t, D_x) = p_m(t, x, D_t, D_x).$

\begin{defn} We say that the Cauchy problem
\begin{equation} \label{eq:1.2}
Pu = f\: {\rm in}\: \Omega \cap \{t < T\},\: {\rm supp}\: u \subset \bar{G}
\end{equation}
is well posed in $G$ if\\
(i) (existence) for every $f \in C_0^{\infty}(\Omega), \: {\rm supp}\: f \subset \overline{\Omega_T^{-}}$ there
exists a solution $u \in {\mathcal E}'(\Omega)$ satisfying $(\ref{eq:1.2})$.\\
(ii) (uniqueness) if $u \in {\mathcal E}'(\Omega)$ satisfies $(\ref{eq:1.2})$, then for every $s, 0 < s \leq T,$ if $Pu = 0$ in $\Omega_s^{-}$, then $u = 0$ in $\Omega_s^{-}$.
\end{defn}
A necessary condition for the well posedeness of the Cauchy problem (WPC) is the hyperbolicity of the operator $P$ at every point  $(t, x) \in G.$
\begin{defn} We say that the operator $P$ with principal symbol $p_m$ is strongly hyperbolic in $G$ if for every point $z_0 = (t_0, x_0) \in G$ there exists a neighborhood $U$ of $z_0$  and $T_0 \geq 0$ ($T_0 < t_0$ if $t_0 > 0$ and $T_0= 0$ if $t_0 = 0$) such that the Cauchy problem $(\ref{eq:1.2})$ for the operator $L = P_m(t, x, D_t, D_x) + Q_{m-1}(t, x, D_t, D_x)$ is well posed in $U^{+}_s$ for every $T_0 \leq s < T(U)$ and for any operator $Q_{m-1}(t, x, D_t, D_x)$ of order less or equal to $m-1$.
\end{defn}
A classical result says that if $P$ is strictly hyperbolic, that is the equation (\ref{eq:1.1}) has simple roots $\lambda_j(t, x,\xi)$ for all $(t. x, \xi)\in G \times \R^N \setminus \{0\}$, then $P$ is strongly hyperbolic. If the equation (\ref{eq:1.1}) has real roots with constant multiplicity for $(t, x, \xi) \in G \times \R^n \setminus \{0\}$, the operator $P$ is strongly hyperbolic {\bf if and only if} it is strictly hyperbolic.  Thus if we have some roots with constant multiplicity $m_j \geq 2$ for the (WPC) we must impose some conditions on lower terms $Q_{m-1}$ called Levi conditions. The analysis of the Cauchy problem for such operators is complete and we know the necessary \cite{FS} and sufficient \cite{Ch} conditions for (WPC).\\

Passing to the case of variable multiplicity of the roots of (\ref{eq:1.1}), notice that the roots $\lambda_j(t, x, \xi)$ in general are not smooth but only continuous. The case of operators with constant coefficients is also completely examined and $P$ is strongly hyperbolic {\bf if and only if} $P$ is strictly hyperbolic. The necessary and sufficient condition of G\aa rding \:for (WPC) says that there exists a constant $c > 0$ such that for the full symbol $p$ of $P$ we have
$$p(\tau, \xi) \neq 0,\:{\rm for}\: |{\rm Im}\: \tau| > c,\: \forall \xi \in \R^n.$$

To understand the situation of variable multiplicity and variable coefficients, consider the example
\begin{equation}
P = D_t^2 - a(z) D_x^2 + b_0(z) D_t + b_1(z)D_x + c(z),\: z = (t,x) \in \R^2
\end{equation}
with $a(z) \geq 0.$ If $a(z_0) = d a(z_0) = a_{tt}(z_0) = 0,\: b_1(z_0) \neq 0,$ in a point $z_0 \in G,$ the Cauchy problem for $P$ is not well posed. On the other hand, if for a point $z_0 = (t_0, x_0) \in G,$ we have $a(z_0) = d a (z_0) = 0,\: a_{tt}(z_0) \neq 0$, then there exists a neighborhood $U$ of $z_0$ such that   the Cauchy problem in $U^{+}_{t_0}$
is well posed for arbitrary smooth lower order terms \cite{O} and $u \in H^{k + 2}(U)$ if $f \in H^{k + N}(\R^2),\: k \in N,$
where 
$$N = 3 + 2 \Bigl[ \frac{3}{2} + \Bigl| b_1(z_0) \Bigr(a_{tt}(z_0)\Bigr)^{-1/2}\Bigr| \Bigr],$$
$[z]$ being the integer part of $z$.\\

 Below we change the notations and we denote $t = x_0,\: x = (x_0, x_1,...,x_n) \in \R^{n+1}.$ The dual variables will be denoted by $\xi = (\xi_0, \xi_1,...,\xi_n) =(\xi_0, \xi').$ Let $\Sigma(p) = \{z \in T^*{\Omega} \setminus \{0\}:\: p(z) = 0\},\: \Sigma_1(p) = \{ z \in T^*(\Omega): z \in \Sigma(p),\:\: dp(z) = 0\}.$ If we have a critical point $(\hh, \hhx) \in \Sigma_1(p)$, then the Hamiltonian system
$$\frac{dx}{ds} = \partial_{\xi} p,\: \frac{d \xi}{ds} = -\partial_{x} p$$
has a stationary point and we consider the differential of the right hand part. Thus we obtain the fundamental matrix
$$F_{p} (\hat{x}, \hat{\xi}) = \Bigl(\begin{matrix} p_{\xi, x}(\hat{x}, \hat{\xi}) & &  p_{\xi, \xi}(\hat{x}, \hat{\xi})\\
-p_{x, x}(\hat{x}, \hat{\xi}) & &  - p_{x, \xi}(\hat{x},\hat{\xi})\end{matrix}\Bigr).$$

We note below two properties of $F_p$:\\
1. For every  point $z \in\Sigma_1(p)$ the Hessian $Q_p(X, Y),\: X, Y \in T_z(T^*(\Omega))$ at $z$ of $\frac{p}{2}$ is well defined. Then
$ Q_p(X, Y) = \sigma(X, F_p(z) Y),$ $\sigma$ being the symplectic form on $T^*(\Omega).$
Thus after canonical transformation the fundamental matrix is transformed into a similar one and its eigenvalues are invariant under canonical transformations. H\"ormander \cite{H1} called $F_p(z)$ Hamiltonian map of $Q_p$.\\
2. If $P$ is hyperbolic in $G$ and $(\hat{x}, \hat{\xi})$ is a critical point of $p_m(x, \xi)$, then $F_{p_m}(\hat{x}, \hat{\xi})$ has at most two non vanishing real simple eigenvalues $\mu$ and $-\mu$ and all other eigenvalues $\lambda$ are purely imaginary, that is ${\rm Re}\: \lambda = 0.$\\
 
The existence of non vanishing real eigenvalues of $F_{p_m}(\hh, \hhx)$ is a necessary condition for strong hyperbolicity. More precisely, let $p_{m-1}(x, \xi) = \sum_{|\alpha| = m-1} c_{\alpha}(x) \xi^{\alpha}$ and let 
$$p_{m-1}'(x, \xi) = p_{m-1}(x, \xi) + \frac{\bf i}{2} \sum_{j=0}^n \frac{\partial^2 p_{m}}{\partial x_j \partial\xi_j} (x, \xi)$$
 be the subprincipal symbol of $P$ which is invariantly defined for $(x, \xi) \in \Sigma_1(p_m).$ Then we have the following

\begin{thm} [\cite{IP}] If $P$ is strongly hyperbolic in $G$, then at every  point $(\hh, \hhx) \in \Sigma_1(p_m)$ the fundamental matrix $F_{p_m}(\hat{x}, \hat{\xi})$ has two non-zero real eigenvalues. Moreover, for $(x, \xi') \in \overset{\circ}G \times (\R^n \setminus \{0\})$ the multiplicities of the roots of (1) are not greater than two, and for $(x, \xi') \in \{x_0 = 0\} \times \R^n \setminus \{0\}$ or for $(x, \xi') \in \{x_0 = T\} \times \R^n \setminus \{0\}$ these multiplicities are not greater than three.
If $F_{p_m}(\hh, \hhx)$ has only purely imaginary eigenvalues, the condition ${\rm Im}\: p_{m-1}'(\hh, \hhx) = 0$  is necessary for (WPC).

\end{thm}                                       
If $F_{p_m}(\hh, \hhx)$ has only purely imaginary eigenvalues, for (WCP) we have a second necessary
condition 
$$|{\rm Re}\: p_{m-1}'(\hh, \hhx)| \leq \frac{1}{4} \sum_{j=0}^{2n + 2}|\mu_j|,$$
$\mu_j$ being the eigenvalues of $F_{p_m}(\hh, \hhx)$ repeated following their multiplicities. This condition has been proved in \cite{IP} in some special cases concerning the structure of $F_{p_m}(\hh, \hhx)$ and without any restriction by H\"ormander \cite{H1}.
\begin{defn} A hyperbolic operator with principal symbol $p(x, \xi)$ will be called effectively hyperbolic if at every point $(\hh, \hhx) \in \Sigma_1(p)$, the fundamental matrix $F_p(\hh, \hhx)$ has two non-zero real eigenvalues.
\end{defn} 
V. Ivrii introduced the following\\

{\bf Conjecture} {\it A hyperbolic operator is strongly hyperbolic if and only if it is effectively hyperbolic.}\\

For operators with at most double characteristics some results for special class of operators have been obtained by H\"ormander \cite{H1}, Ivrii \cite{I} and Melrose \cite{M}. The sufficient part of the above conjecture is difficult since the double roots of the equation (1) in general are not smooth and we have not a factorization with smooth factors. Moreover, the loss of regularity could depend on the point and a microlocalization leads to considerable difficulties when we must treat the commutators. The above conjecture for operators with double characteristics has been completely solved  by N. Iwasaki \cite{Iw1}, \cite{Iw2} and T. Nishitani \cite{N1}, \cite{N2}.
The proofs are rather long and very technical.\\

 An effectively  hyperbolic operator could be strongly hyperbolic if it has triple characteristics on the boundary on $G$ but to our best knowledge there are no examples of such operators in the literature. Our purpose is to study a class of operators $P$ with triple characteristics on $t = 0$ and to prove that $P$ is strongly hyperbolic. Thus the above conjecture is true for some special operators with triple characteristics. The analysis of the general case remains open. 

\section{Hyperbolic operators with triple characteristics}
 In this section we use again the notations of Section 1. According to Theorem 1,
an effectively hyperbolic operator $P$ in $G$ may have triple characteristics in $G$ only for $t = 0$ or $t = T$. Assume that $P$ has triple characteristics for $t = 0$ and suppose that the triple roots of (\ref{eq:1.1}) for $t = 0$ are $\tau = 0$ (in general the triple characteristics for $t = 0$ are $ \tau = \lambda(0, x,\xi)$). Let $P$ be of order 3 and let 
$$p_3 = \tau^3 + q_1(t, x, \xi)\tau^2 + q_2(t, x, \xi)\tau + q_3(t, x, \xi)$$
be the principal symbol of $P$ with $q_j,\: j = 1, 2, 3,$ real-valued polynomials of order $j$ with respect to $\xi$ with smooth coefficients.

\begin{lem}[\cite{IP}] Let $p_3(t, x, \tau, \xi)$ be hyperbolic in $G$ and let $\tau = 0$ be a triple root of $p_3(0, x, \tau, \xi) = 0,\:(0, x) \in G$. Then
$$q_3(0, x, \xi) = \partial_t q_3(0, x, \xi) = q_2(0, x, \xi) = q_1(0, x, \xi) = 0,\: (0,x) \in G, \: \xi \in \R^n.$$
Moreover, $p_3$ is effectively hyperbolic  for $t = \tau = 0$, if and only if
$$\frac{\partial^2 p_{3}}{\partial\tau \partial t} (0, x, 0, \xi) < 0,\: \xi \in \R^n \setminus \{0\}.$$
\end{lem}

Thus we must study an operator $P$ with principal part
$$P_3 = D_t^3 +  t a_1(t, x, D_x)D_t^2 - t a_2(t, x, D_x)D_t + t^2 a_3(t, x, D_x)$$
with $a_j(t, x, \xi)$ real-valued polynomials of order $j$ in $\xi$ and  $a_2(t, x, \xi) \geq c|\xi|^2,\: c > 0$ for $\xi \neq 0.$ We write $P = P_3 + Q$
with lower order terms $Q = B_2(t,x, D_x) + B_1(t, x, D_x) D_t + C(t, x, D_t, D_x).$ Here $B_2$ and $B_1$ are differential operator of order 2 and 1, respectively, while $C$ is an operator of order 1.
Notice that for $|\xi| = 1$ the discriminant $\Delta$ of the equation $p_3(t, x, \tau, \xi) = 0$ with respect to $\tau$ has the form
$$\Delta(t, x, \xi) =  \Bigl(\frac{-3ta_2 - t^2a_1^2}{9}\Bigr)^3 + \Bigl(\frac{-9t^2a_1 a_2- 27 t^2 a_3 - 2 t^3a_1^3}{54}\Bigr)^2$$
$$ = q^3 + r^2 = -\frac{1}{27} t^3 a_2^3 + {\mathcal O}(t^4)a_6$$
and $\Delta \leq 0$ for small $t \geq 0.$ Thus the operator $P$ is strictly hyperbolic for small $t > 0$ and it suffices to examine the Cauchy problem for $0 \leq t \leq t_0,\: t_0 \ll 1$. Since the coefficients of the cubic equation $p_3(t,x,\tau, \xi) = 0$  are real, for $t \geq 0$ its real roots $\lambda_k(t, x,\xi),\: k =1,2,3,$ have the following trigonometric form (see for instance, \cite{ SL})
$$\begin{cases} \lambda_1 = 2\rho^{1/3}\cos(\theta/3) -\frac{t a_1}{3},\\
\lambda_2 = 2\rho^{1/3} \cos(\theta/3 + \frac{2\pi}{3}) - \frac{t a_1}{3},\\
\lambda_3 =  2\rho^{1/3}\cos(\theta/3 + \frac{4 \pi}{3}) - \frac{t a_1}{3},\end{cases}$$
where 
$$\rho = ( -q)^{3/2},\: \theta = \arccos (r/\rho).$$
Next consider the symbols
$$\delta_k = \frac{\partial p_3}{\partial \tau} \Bigl\vert_{\tau = \lambda_k} = \Bigl(3\tau^2 + 2ta_1 \tau - t a_2\Bigr)\Bigl\vert_{\tau = \lambda_k},\: k = 1, 2, 3.$$
Since these symbols are homogeneous of order 2 in $\xi$, to find lower bounds for $|\delta_k|$, it is sufficient to examine their behavior for $|\xi| = 1.$
We have 
$$\delta_1 = 12 \rho^{2/3} \cos^2(\theta/3) - ta_2 + {\mathcal O}(t^{3/2})a_2= \Bigl(4 \cos^2(\theta/3) - 1\Bigr) t a_2+ {\mathcal O}(t^{3/2})a_2.$$
Since $\frac{r}{\rho} = {\mathcal O}(t^{1/2})$, we have  $\cos(\theta/3) =\frac{\sqrt{3}}{2} + o(t)$ and this implies
for small $t$ and $|\xi| = 1$ the estimate $|\delta_1| \geq c_1 t a_2$ with $c_1 > 0.$ On the other hand,

$$\delta_{2,3} = 3 \lambda_{2,3}^2 - ta_2 + {\mathcal O}(t^{3/2}) a_2 = \Bigl(4\sin^2(\pi/6 \pm \theta/3) - 1\Bigr)ta_2 + {\mathcal O}(t^{3/2})a_2$$
and we obtain the following
\begin{lem} There exist constants $\gamma > 0$ and $\gamma_1 > 0$ such that for $0 \leq t \leq \gamma_1$ we have
\begin{equation}
|\delta_k| \geq \gamma\: t a_2(t,x,\xi) \geq \gamma c t |\xi|^2,\: k = 1,2,3.
\end{equation}
\end{lem} 
Finally, notice that $\lambda_1 \lambda_2 \lambda_3 = - t^2 a_3(t, x, \xi).$

%Consider a scaling $t = \epsilon^{2/3}s,\: x = \epsilon y,\: \epsilon > 0.$ Multiplying by $\epsilon^2$, we obtain %an operator
%$$P = D_s^3 - s a_2(\epsilon^{2/3} s, \epsilon y, D_y)D_s^2 + B_2(\epsilon^{2/3}s, \epsilon y, D_y) $$
%$$+ \epsilon^{1/3} \Bigl[s  a_1(\epsilon^{2/3} t, \epsilon y, D_y)D_s^2 + s^2a_3(\epsilon^{2/3}s, \epsilon y, D_y) + 
%B_1(\epsilon^{2/3} s, \epsilon y, D_y)D_s\Bigr] + \epsilon C_1(...).$$
%We choice $\epsilon = {\mathcal O}(\frac{1}{N})$, where $N$ is a big fixed integer related to 
%$$1 + \sup_{0 \leq t \leq t_0,x \in U, |\xi| = 1}  \Bigl[| B_2(t, x, \xi) a_2(t, x,\xi)^{-1}|\Bigr].$$
%We will study a model operator keeping the leading terms in the above representation.

\section{Energy estimates for a model operator}
Consider the operator
\begin{equation}
\label{1}
\def\aa{a_{2}}
P(t, D_{t}, D_{x}) = D_{t}^{3} + ta_1(t,D_x)D_t^2- t\aa(D_x)D_t  + t^2 a_3(t, D_x) + b(t, D_x), \: t \geq 0
\end{equation}
where $a_2(D_x) = \sum_{i,j=1}^n a_{i,j} D_iD_j$ and  $b(t, D_x) = \sum_{i,j=1}^{n} b_{i,j}(t) D_iD_j$ is a second order differential operator. For simplicity we assume that $a_2$ is independent on $t$. The analysis of operators with $a_2(t,D_x)$ goes without any change. We assume that 
$$a(\xi) = a_2(\xi) = \sum_{i,j =1}^n a_{i,j} \xi_i \xi_j \geq \delta_0 |\xi|^2,\: \delta_0 > 0.$$
Moreover, the symbols $a_1(t, \xi),\: a_3(t, \xi)$ are real-valued and homogeneous of order 1 and 3 in $\xi$, respectively.
 We want to
establish an \textit{a priori} estimate for $ P $ for $ t \geq 0 $. Set 
$$f(t, \xi) = t + \frac{1}{(1 + a(\xi))^{1/3}}.$$
Let $v(t, x) \in C_0^{\infty}(\R_t \times \R^n)$. Multiplying $P $ by $ -\ii $ and taking the Fourier transform
with respect to the variable $ x $, we obtain 
$$ \hat{P}u = \widehat{-\ii Pv} = \partial^{3}_{t}u + \ii t a_1(t, \xi) u'' +
t \partial_{t} a(\xi)u  -\ii t^2 a_3(t, \xi) u + \bv$$
with $b_1(t,D_x) = -\ii b(t,D_x)$ and $u = \hat{v}$. Let $u'' = \widehat{v_{tt}},\: u' = \hat{v_t}.$
We have
$$ 
2 \re \hat{P} u \bar{u}'' = \partial_{t} |u''|^{2} + t a(\xi) \partial_{t}
|u'|^{2} + 2 \ia %(t^2 a_3)\im (u \bar{u}'')% 
+ 2 \re \Bigl(\bv\bar{u}''\Bigr).
$$
Denote by $N$ a large positive integer and by $ \lambda $ a large
positive parameter.  Multiply the above
identity involving $ \hat{P}u $ by $ e^{-\lambda t} f^{- 2N}
$. We obtain
\begin{multline*}
e^{-\lambda t} f^{-2N} 2 \re(\hat{P}u \bar{u}'') = e^{-\lambda t}
f^{-2N} \partial_{t} |u''|^{2} + e^{-\lambda t}f^{-2N} t
a(\xi) \partial_{t} |u'|^{2} 
\\
+ e^{-\lambda t} f^{-2N} 2 \Bigl( \ia + \re \bv \bar{u}'' \Bigr)
\\
= e^{-\lambda t} f^{-2N} \partial_{t} \tilde{E}(u) - e^{-\lambda t}
f^{-2N} a(\xi) |u'|^{2} + e^{-\lambda t} f^{-2N} 2 \Bigl(\ia + \re \bv \bar{u}''\Bigr),
\end{multline*}
where 
$$ 
\tilde{E}(u) = |u''|^{2} + t a(\xi) |u'|^{2}.
$$
The above identity can be rewritten as
\begin{multline*}
e^{-\lambda t} f^{-2N} 2 \re(\hat{P}u \bar{u}'') = 
\partial_{t} \left( e^{-\lambda t} f^{-2N} \tilde{E}(u) \right) +
\lambda e^{- \lambda t} f^{-2N} \tilde{E}(u) 
\\
+ 2N e^{-\lambda t} f^{- 2N -1} \tilde{E}(u) - e^{-\lambda t} f^{-2N}
a(\xi) |u'|^{2} + 2e^{-\lambda t} f^{-2N}  \Bigr( \ia + \re \bv \bar{u}''\Bigr).
\end{multline*}
Since
$$ 
e^{-\lambda t} f^{2N} 2 \re(\hat{P}u \bar{u}'') \leq e^{-\lambda t}
f^{-2N+1} |\hat{P}u |^{2} + e^{-\lambda t} f^{-2N-1} | u''|^{2},
$$
we have the inequality
\begin{multline*}
e^{-\lambda t} f^{-2N+1} | \hat{P}u |^{2} \geq \partial_{t}\left(
  e^{-\lambda t} f^{-2N} \tilde{E}(u) \right) 
+ \lambda e^{-\lambda t} f^{-2N} \tilde{E}(u) 
\\
+ (2N-1) e^{-\lambda t}
f^{-2N-1} |u''|^{2} + 2N e^{-\lambda t} f^{-2N-1} t a(\xi) |u'|^{2} 
\\
- e^{-\lambda t} f^{-2N} a(\xi) |u'|^{2} + 2e^{-\lambda t} f^{-2N} \Bigr( \ia + 
\re \bv \bar{u}''\Bigr).
\end{multline*}

Let us now consider the following identity, where $ k $ is a positive
integer and $g$ denotes a smooth function in the same class as $ u
$:

$$e^{-\lambda t} f^{-2k} 2 \re g' \bar{g} = \partial_{t}
\left(e^{-\lambda t} f^{-2k} |g|^{2} \right) + \lambda e^{-\lambda t}
f^{-2k} |g|^{2} + 2 k e^{-\lambda t} f^{-2k-1} |g|^{2}.$$
This implies
$$e^{-\lambda t} f^{-2k+1} |g'|^{2} \geq  \partial_{t}
\left(e^{-\lambda t} f^{-2k} |g|^{2} \right) + \lambda e^{-\lambda t}
f^{-2k} |g|^{2} + (2 k - 1)  e^{-\lambda t} f^{-2k-1} |g|^{2}.$$

Now, taking $g = u'$ we have

\begin{equation}
\label{3}
e^{-\lambda t} f^{-2k+1} |u''|^{2} \geq  \partial_{t}
\left(e^{-\lambda t} f^{-2k} |u'|^{2} \right) + \lambda e^{-\lambda t}
f^{-2k} |u'|^{2}
+ (2 k - 1)  e^{-\lambda t} f^{-2k-1} |u'|^{2},
\end{equation}
while, taking $g = u $, we get
\begin{equation}
\label{4}
e^{-\lambda t} f^{-2k+1} |u'|^{2} \geq  \partial_{t}
\left(e^{-\lambda t} f^{-2k} |u|^{2} \right) + \lambda e^{-\lambda t}
f^{-2k} |u|^{2} + (2 k - 1)  e^{-\lambda t} f^{-2k-1} |u|^{2}.
\end{equation}
From (\ref{3}) and (\ref{4}) above we obtain
\begin{multline}
\label{5}
e^{-\lambda t} f^{-2k+1} |u''|^{2} \geq \partial_{t} \left(
  e^{-\lambda t} f^{-2k} |u'|^{2}\right) + \lambda e^{-\lambda t}
f^{-2k} |u'|^{2} 
\\
+ (2k-2) e^{-\lambda t} f^{-2k-1} |u'|^{2} 
\\
+ \partial_{t} \left(e^{-\lambda t} f^{-2k-2} |u|^{2}\right) + \lambda
e^{- \lambda t} f^{-2k-2} |u|^{2} + (2k+1) e^{-\lambda t} f^{-2k-3}
|u|^{2}. 
\end{multline}
Plugging this into the estimate for $ |\hat{P}u |^{2} $ and choosing $
k = N+1 $, we obtain
\begin{multline}
\label{6}
e^{-\lambda t} f^{-2N+1} |\hat{P} u|^{2} \geq \partial_{t} \left(
  e^{-\lambda t} f^{-2N} \tilde{E}(u) \right) + \lambda e^{-\lambda t}
f^{-2N} \tilde{E}(u)
\\
+ {\mathcal O}(N) \left \{ e^{-\lambda t} f^{-2N-1} |u''|^{2}
  + \partial_{t} \left( e^{-\lambda t} f^{-2N-2} |u'|^{2} \right) +
  \lambda e^{-\lambda t} f^{-2N-2} |u'|^{2} \right\} 
\\
+ {\mathcal O}(N^{2}) e^{-\lambda t} f^{-2N-3} |u'|^{2} 
\\
+ {\mathcal O}(N) \left \{ \partial_{t} \left(e^{-\lambda t} f^{-2N-4}
    |u|^{2}\right) + \lambda e^{-\lambda t} f^{-2N-4} |u|^{2} \right\}
\\
+ {\mathcal O}(N^{2}) e^{-\lambda t} f^{-2N-5} |u|^{2}
\\
+ 2N e^{-\lambda t} f^{-2N-1} t a(\xi) |u'|^{2} - e^{-\lambda t}
f^{-2N} a(\xi) |u'|^{2}
\\
+ 2e^{-\lambda t} f^{-2N} \Bigr( \ia + \re \bv \bar{u}''\Bigr).
\end{multline}
Here $ {\mathcal O}(N) $ means a function of $ N $ which satisfies an
estimate of the type: $ {\mathcal O}(N) \geq c N $, with a
\textit{fixed} positive constant $ c $.

From inequality (\ref{5}) above we also deduce that
\begin{multline*}
e^{-\lambda t} f^{-2N-1} t a(\xi) |u'|^{2} \geq \partial_{t}\left(
  e^{-\lambda t} f^{-2N-2} t a(\xi)|u|^{2}\right) + \lambda e^{-\lambda t}
f^{-2N-2} t a(\xi) |u|^{2} 
\\
- e^{-\lambda t} f^{-2N-2} a(\xi) |u|^{2} + (2N+1) e^{-\lambda t}
f^{-2N-3} t a(\xi) |u|^{2}.
\end{multline*}
Replacing the part of the corresponding term in (\ref{6}) with the
above inequality, we finally obtain
\begin{multline}
\label{7}
\el f^{-2N+1} |\hat{P} u|^{2} \geq \partial_{t} \left(
  e^{-\lambda t} f^{-2N} \tilde{E}(u) \right) + \lambda e^{-\lambda t}
f^{-2N} \tilde{E}(u)
\\
+ {\mathcal O}(N) \left \{ e^{-\lambda t} f^{-2N-1} |u''|^{2}
  + \partial_{t} \left( e^{-\lambda t} f^{-2N-2} |u'|^{2} \right) +
  \lambda e^{-\lambda t} f^{-2N-2} |u'|^{2} \right\} 
\\
+ {\mathcal O}(N^{2}) e^{-\lambda t} f^{-2N-3} |u'|^{2} 
\\
+ {\mathcal O}(N) \left \{ \partial_{t} \left(e^{-\lambda t} f^{-2N-4}
    |u|^{2}\right) + \lambda e^{-\lambda t} f^{-2N-4} |u|^{2} \right\}
\\
+ {\mathcal O}(N^{2}) e^{-\lambda t} f^{-2N-5} |u|^{2}
+ {\mathcal O}(N) e^{-\lambda t} f^{-2N-1} t a(\xi) |u'|^{2} 
\\
+ {\mathcal O}(N) \left\{ \partial_{t}\left( e^{-\lambda t} f^{-2N-2}t a(\xi)
    |u|^{2}\right) + \lambda e^{-\lambda t} f^{-2N-2} t a(\xi)
  |u|^{2}\right\}
\\
+ {\mathcal O}(N^{2}) e^{-\lambda t} f^{-2N-3} t a(\xi) |u|^{2}  
\\
- {\mathcal O}(N) e^{-\lambda t} f^{-2N-2} a(\xi) |u|^{2}
- e^{-\lambda t} f^{-2N} a(\xi) |u'|^{2}\\
+ 2e^{-\lambda t} f^{-2N} \Bigl( \ia + \re \bv \bar{u}''\Bigr).
\end{multline}
There are four ''error'' terms, all written in the last two lines of
(\ref{7}). We deal first with the term containing $ u' $, the second term
term in the second line from below.
Neglecting the exponential term, we would like to
estimate $ f^{-2N} a(\xi) $ by $f^{-2N-3} + f^{-2N-1} t a(\xi)
$. First we would like to prove an inequality of the form
\begin{equation} \label{8}
 \frac{f^{-2N -3}}{1 + a(\xi)} + t f^{-2N-1} \geq \alpha f^{-2N},
\end{equation}
with a positive constant $ \alpha $. Dividing by $ f^{-2N-3} $, the proof is reduced to the inequality 
$$ 
\frac{1}{1 + a(\xi)} + t f^{2} \geq \alpha f^{3}.
$$
Now 
$$ f^{3} = t^{3} + \frac{1}{1 + a(\xi)} + \frac{3t^{2}}{(1 + a(\xi))^{1/3}} + \frac{3t}{(1 + a(\xi))^{2/3}},$$
 while on the left hand side we have 

$$ \frac{1}{1 + a(\xi))} + t^{3} + \frac{2 t^{2}}{(1 + a(\xi))^{1/3}}
 + \frac{t}{(1 + a(\xi))^{2/3}}. $$

 The terms on both sides are the same, so that
if we choose $ \alpha $ suitably, (\ref{8}) ensues. Thus we deduce
\begin{equation} \label{9}
\alpha f^{-2N} a(\xi) \leq ta(\xi) f^{-2N -1} + \frac{a(\xi)}{(1 + a(\xi)} f^{-2N - 3}
\leq ta(\xi) f^{-2N -1} + f^{-2N-3}.
\end{equation}
%Note that we did use only a fixed and possibly small portion of the
%large coefficients of the terms containing $ f^{-2N-3} $ and $
%f^{-2N-1} t a(\xi)$. 

Next let us treat the first term in the second line from below in (\ref{7}). We want to estimate $ f^{-2N-2} a(\xi) $ with $ f^{-2N-5} + f^{-2N-3}
t a(\xi)$. This is very easy, since the coefficients of the terms
containing $ |u|^{2} $ in (\ref{7}) grow as $ N^{2} $, and a small
portion of them may absorb ${\mathcal O}(N) $. Now the inequality
$$ f^{-2N-5} + f^{-2N-3} t a(\xi) \geq \alpha f^{-2N-2} a(\xi)$$
is obtained from (\ref{9}), dividing by $ f^{2} $. 
%In the same way we deal with the second term in the last line of (\ref{7}) estimating
 %$f^{-2N} a(\xi)$ by $f^{-2N -1} t a(\xi) + f^{-2N - 3}$. \\

Now we pass to the analysis  of the last term in the last line of (\ref{7}).  First we deduce
$$\re (\bv \buu) = \re b_1(t,\xi) \re (u\buu) - \im b_1(t,\xi)\im (u \buu).$$
To deal with the term involving $\re b_1(t, \xi)$, we use the equality
$$2 \re (u \buu) = \partial_t 2 \re(u \bu) - 2 |u'|^2$$
The term with $|u'|^2$ be be treated as above since $|b_1(t, \xi)| \leq C\delta_0 a(\xi)$. To study the term with $\re(u \bu)$, we write
\begin{eqnarray} \label{eq:3.9}
\el f^{-2N}\re \Bigl(b_1(t, \xi) \re (u \bu)\Bigr)\\ \nonumber
= \partial_t\Bigl(\el f^{-2N} \re b_1(t, \xi) \re (u \bu)\Bigr) + \lambda \el f^{-2N} \re b_1(t, \xi) \re (u \bu)\\   \nonumber
+ 2N f^{-2N-1} \re b_1(t, \xi) \re(u \bu) + \el f^{-2N} \re b_{1, t}(t, \xi) \re (u\bu) = \partial_t(...) + I + II + III.
\end{eqnarray}
There are three terms on the right hand side of (\ref{eq:3.9}). Consider $I$. Applying the Cauchy-Schwartz
inequality and $|b_1(t, \xi)| \leq C|\xi|^2$, we obtain

\begin{eqnarray}
\lambda\Big|\el  f^{-2N} \re b_1(t, \xi) \re(u \bu)\Big|\\ \nonumber
\leq C \lambda\delta_0^{-1}\Bigl[\epsilon \el f^{-2N + 1}a(\xi)|u'|^2  + \frac{1}{\epsilon} \el   f^{-2N- 1} a(\xi) |u|^2\Bigr]
\\ \nonumber
\leq C \lambda\delta_0^{-1} \epsilon \alpha^{-1}\el [f^{-2N}t a(\xi)|u'|^2 + f^{-2N -2}|u'|^2]\\ \nonumber
+ \frac{C \lambda}{\epsilon}\delta_0^{-1}\alpha^{-1}\el [ f^{-2N-2} t a(\xi) |u|^2 + f^{-2N -4}|u|^2],
\end{eqnarray}
where $ \epsilon > 0$ is a small positive constant, to be chosen below. Taking $C \delta_0^{-1}\alpha^{-1} \epsilon < 1/2$, we may estimate the term with $f^{-2N} t a(\xi)|u'|^2$ by $f^{-2N} \tilde{E}(u)$. Next 
$\frac{C \lambda}{\epsilon}\delta_0^{-1}\alpha^{-1}\el f^{-2N-2} t a(\xi) |u|^2$  can be absorbed by the corresponding term in (6) with large $N$ and the same is true for the term with $f^{-2N-4}|u|^2.$ The analysis of $III$ is similar and simpler.\\

To handle $II$, we use the inequality
$$II \leq C_1^2\delta^{-1} \el  f^{-2N} a(\xi) |u'|^2 + 4N^2 \delta f^{-2N - 2} a(\xi) |u|^2,$$
where $C_1 = C \delta_0^{-1}$ and $\beta > 0$ is a small constant.

The latter term in the above line is similar to the first in the last
line of (\ref{7}); the only difference is the factor in front, which
is bigger here. 
However, remarking that all the terms containing $ |u|^{2} $ in
(\ref{7}) have also ${\mathcal O}(N^{2}) $, it is clear that choosing
$ \delta $ suitably small, but finite and independent of $ u $, $ N $
and $ \lambda $, will allow us to conclude by arguing as above. The
fist summand on the other hand is similar to the middle term in the
last line of (\ref{7}): $ C_1 $ is real and depends on the lower
order terms, $ \delta $ is fixed. This is estimated as we did before,
provided that $ N $ is large enough. \\

Next we turn to the term containing  $ - \im b_1(t, \xi) \im (u
\bar{u}'') $ containing $\im b_1(t, \xi).$ We remark that $ \im (u \bar{u}'') = \partial_{t} (u
\bar{u}' - u' \bar{u}) $, so that we obtain two terms which can be
discussed almost verbatim as before. This might require enlarging $ N
$.\\

Finally, consider the term 
$$2\el f^{-2N} \ia \geq -C_1\el f^{-2N} t^4 (1+|\xi|^2)^3 |u|^2 -  \el f^{-2N} |u''|^2.$$
The last term in the right hand side can be treated as above, however the first one cannot be absorbed by other positive terms taking $N$ large enough. Consequently, in (\ref{7}) we obtain an upper bound on the left by
$$\el f^{-2N+1} |\hat{P} u|^{2} + C_1 t^4\el f^{-2N} (1 + |\xi|^2)^3 |u|^2.$$

Now assume $0 \leq s < T \leq 1$ and $v = D_t v = D_t^2 v = 0$ when $t = s$. We integrate from $t$ to $T $ w.r.t. to the time variable $s$. Thus yields some integrals $\int_t^T (....) ds$ and terms
$$\eT f^{-2N}(T,\xi) \tilde{E}(u(T, \xi))$$
$$ + {\mathcal O}(N)\Bigl[\eT f^{-2N-2}(T, \xi)|u'(T, \xi)|^2 + \eT f^{-2N-4}(T, \xi)|u(T, \xi)|^2\Bigr]$$
$$+ {\mathcal O}(N) \eT f^{-2N-2}(T, \xi) T a(\xi) |u(T, \xi)|^2 $$
$$+ 2\eT f^{-2N}(T, \xi)\re b(T, \xi)\re (u\bu) (T, \xi)$$
$$ - 2\eT f^{-2N} (T, \xi) \im b(T, \xi) (u\bu - \bu u) (T, \xi).$$

Concerning the boundary terms, only the last three terms in the above sum has no positive sign. We have no coefficient
$\lambda$ before them and these terms can be treated by the above argument using the positive terms and choosing $N$ large enough. Next we integrate with respect to $\xi$ in $\R^n$ and we replace the $L^2(\R_{\xi}^n)$ norms by $L^2(\R^n_x)$ norms. Moreover, we have the obvious inequalities
$$f^{-1} = \frac{(1 + a(\xi))^{1/3}}{t(1 + a(\xi))^{1/3} + 1} \leq (1 + a(\xi))^{1/3} \leq C_2(1 + |\xi|^2)^{1/3},\: t \geq 0,$$
$$f^{-1} \geq 1/2, \: 0 \leq t < T \leq 1.$$
Then 
$$C_1 t^4\int_t^T \int e^{-\lambda s} f^{-2N} (1 + |\xi|^2)^3 |u|^2 d\xi ds \leq C_2 t^4\int_t^T\int e^{-\lambda s} (1 + |\xi|^2)^{2N/3 + 3} |u|^2 d\xi ds$$
$$\leq C_3 t^4 \int_t^T e^{-\lambda s} \|v\|_{(2N/3 + 3)}^2 ds,$$
where $\|. \|_{(s)}$ is the $H_{(s)}$ norm in $\R^n$ for fixed $s.$ Now for $v$ and $0 < t \leq T$ and small $T$ we may apply the energy estimates for strictly hyperbolic operators (see Section 23.2 and the proof of Lemma 23.2.1 in \cite{H1}). Taking into account Lemma 2,  we get
$$\int_t^T e^{-\lambda s} \|v\|_{(2N/3 + 3)}^2 ds \leq \frac{C_{N}}{t^2} \int_t^T e^{-\lambda s} \|Pv\|_{(2N/3 + 1)}^2ds.$$

We introduce $U_1(s, \xi) = (1 +|\xi|^2)^{1/2}u(s, \xi)$ and observe that $U_1$ satisfies the same initial conditions on $s = t$ as $u$ and
$$\hat{P}U_1 = (1 + |\xi|^2)^{1/2}\hat{P}u.$$
Finally, we obtain the following
\begin{thm} Let $v \in C_0^{\infty}(\R_t \times \R^n)$ and let $v(s,x) = D_t v(s, x) = D_t^2 v(s, x) = 0$ for $s = t.$
Let $0 \leq t < T \leq 1$ Then for $T$ small enough and for an integer $N$ and $\lambda > \lambda_0$ depending on the lower order terms $b(t, D_x)$ we have the estimate
\begin{equation} \label{eq:3.11}
\lambda \int_t^T e^{-\lambda s} \Bigl(\|D_t^2 v\|_{(1)}^2 + \|D_tv\|_{(2)}^2 + \|v\|_{(2)}^2\Bigr) ds \leq C(N) \int_t^T e^{-\lambda s}\|Pv\|_{(2N/3 +2)}^2 ds.
\end{equation}

\end{thm}

Now will treat the estimates for functions $v \in C_0^{\infty}(\R_t \times \R^n)$ with initial data
$$v(T,x) = D_t v(T, x) = D_t^2 v(T, x) = 0.$$
 To do this we multiply $\hat{P}u$ by $-e^{\lambda t} f^{2N} \bar{u}''$ and 
repeating the above argument, we obtain for $0 \leq t < T \leq 1$
\begin{multline}
\label{6}
 e^{\lambda t} f^{2N+1} |\hat{P} u|^{2} \geq -\partial_{t} \left(
  e^{\lambda t} f^{2N} \tilde{E}(u) \right) + \lambda e^{\lambda t}
f^{2N} \tilde{E}(u)
\\
+ {\mathcal O}(N) \left \{ e^{\lambda t} f^{2N-1} |u''|^{2}
  - \partial_{t} \left( e^{\lambda t} f^{2N-2} |u'|^{2} \right) +
  \lambda e^{\lambda t} f^{2N-2} |u'|^{2} \right\} 
\\
+ {\mathcal O}(N^{2}) e^{\lambda t} f^{2N-3} |u'|^{2} 
\\
+ {\mathcal O}(N) \left \{ -\partial_{t} \left(e^{\lambda t} f^{2N-4}
    |u|^{2}\right) + \lambda e^{\lambda t} f^{2N-4} |u|^{2} \right\}
\\
+ {\mathcal O}(N^{2}) e^{\lambda t} f^{2N-5} |u|^{2}
\\
- 2N e^{\lambda t} f^{2N-1} t a(\xi) |u'|^{2} + e^{\lambda t}
f^{2N} a(\xi) |u'|^{2}
\\
+ 2 e^{\lambda t} f^{2N} \ia - 2 e^{\lambda t} f^{2N}  \re \bv \bar{u}''.
\end{multline}

Now we assume $0 \leq s < T \leq 1$ and let $v = D_tv = D_t^2v = 0$ when $s = T.$ We integrate from $t$ to $T$ with respect to the time variable $s$ and we treat the boundary terms with $s = t$ as above, while the "error" terms
are handled in the same way as in the case with initial data on $s = t.$ Thus we obtain a priori estimate
involving the "weights" $f^{2N -k}(D_x)$, $-1 \leq k \leq 5.$ 
%We replace $u$ by $U_N = (1 +a(\xi))^{-(2N +1)/6} u (t, \xi)$ which satisfy the same boundary conditions on $s = T$ and we observe that $\hat{P}U_N = (1 + a(\xi))^{-(2N + 1)/6} \hat{P}u.$\\
On the other hand,
$$f^{2N +1} \leq (t +1)^{2N + 1},\: 0 \leq t < T \leq 1,$$
$$ f^{2N} \geq \frac{1}{(1 + a(\xi))^{2N/3}} \geq B_N(1 + |\xi|^2)^{-2N/3}. $$
%Estimating the right and the left hand side of the integrals on $[t, T]$, we obtain

%$$\lambda\int_t^T e^{\lambda s}\Bigl((1 + |\xi|^2)^{-2N/3}\tilde{E}(u) + (1 + |\xi|^2)^{-(2N-2)/3}|u'|^2 + (1 + |\xi|^2)^{-(2N-4)/3}|u|^2\Bigr) ds$$
%$$+ \int_t^T e^{\lambda s}(1 + |\xi|^2)^{-2N/3 + 1}|u|^2 ds
%\leq C(N) \int_t^T e^{\lambda s} (s + 1)^{2N + 1}|\hat{P}u|^2 ds.$$

We introduce $U_N(s, \xi) = (1 +|\xi|^2)^{(2N + 2)/3}u(t, \xi)$ and observe that $U_N$ satisfies the same initial conditions on $s = T$ as $u$ and
$$\hat{P}U_N = (1 + |\xi|^2)^{(2N + 2)/3}\hat{P}u.$$

%Applying the above estimate for $U_N$, we get
%$$\lambda\int_t^T e^{\lambda s}\Bigl(\tilde{E}(u) + (1 + |\xi|^2)^{4/3}|u'|^2 + (1 + |\xi|^2)^2|u|^2\Bigr) ds$$
%$$\leq C(N) \int_t^T e^{\lambda s} (s + 1)^{2N + 1}(1 + |\xi|^2)^{(2N +2)/3}|\hat{P}u|^2 ds,$$
Thus we deduce the following
\begin{thm}Let $v \in C_0^{\infty}(\R_t \times \R^n)$ and let $v(s,x) = D_t v(s, x) = D_t^2 v(s, x) = 0$ for $s = T.$
Let $0 \leq t < T \leq 1.$ Then for $T$ small enough and for an integer $N$ and $\lambda > \lambda_0$ depending on the lower order terms $b(t, D_x)$ we have the estimate
\begin{equation}
\lambda \int_t^T e^{\lambda s} \Bigl(\|D_t^2 v\|_{(2/3)}^2 + \|D_tv\|_{(4/3)}^2 + \|v\|_{(2)}^2\Bigr) ds
 \leq C_1(N, T) \int_t^T e^{\lambda s}\|Pv\|_{(2N/3 + 2)}^2 ds,
\end{equation}
where $\|. \|_{(m)}$ is the $H_{(m)}$ norm in $\R^n$ for fixed $s.$
\end{thm}

From Theorems 2 and 3 we conclude in a standard way that the Cauchy problem for $P$ is well posed.

\section{Operators with coeffcients depending on $t$ and $x$}

We sketch briefly some ideas for the analysis of the case when we have operators with coefficients depending on $t$ and $x$. 

First consider a scaling $t = \epsilon^{2/3}s,\: x = \epsilon y,\: \epsilon > 0.$ Multiplying by $\epsilon^2$, we obtain an operator
$$P = D_s^3 - s a_2(\epsilon^{2/3} s, \epsilon y, D_y)D_s + B_2(\epsilon^{2/3}s, \epsilon y, D_y) $$
$$+ \epsilon^{1/3} \Bigl[s  a_1(\epsilon^{2/3} t, \epsilon y, D_y)D_s^2 + s^2a_3(\epsilon^{2/3}s, \epsilon y, D_y) + 
B_1(\epsilon^{2/3} s, \epsilon y, D_y)D_s\Bigr] + \epsilon C_1(...).$$

Our final purpose is to choose $\epsilon = {\mathcal O}(\frac{1}{N})$, where $N$ is a big fixed integer related to 
lower order terms as in the case treated in Section 3. With this
choice of $\epsilon$ we are going to study the Cauchy problem for
sufficiently small $t > 0$. This is enough since for $t > 0$ our operator is strictly hyperbolic.

\bigskip

We cannot apply Fourier transform and moreover it is convenient to employ a suitable class of pseudodifferential operators. Notice that $ f = t + (1 + a_2(t,x,\xi))^{-1/3}$ is a symbol in the class $S^{0}_{1,2/3}$, when derivatives with respect to $t$ are considered, but is in the class $S^{0}_{1, 0}$ if $ t $ is just a parameter and no
derivatives with respect to $t$ are involved.

Let $ \langle \xi \rangle^{2} = 1 + |\xi|^{2}$ and let
\begin{equation}\label{eq:10metric}
g_{(x, \xi)} = |dx|^{2} + \langle \xi \rangle^{-2} |d\xi|^{2}
\end{equation}
be the classical slowly varying $(1, 0)-$ metric. We need also the
dilated metric 
\begin{equation}
\label{eq:dilatedmetric}
g_{(x, \xi)}^{\epsilon} = \epsilon^{2} |dx|^{2} + \langle \xi \rangle^{-2} |d\xi|^{2}.
\end{equation}

Define the following ''order'' function
\begin{equation}
\label{eq:4.1}
m_{N}^{t, \mu}(x, \xi) = f^{-N}(t, \xi) \langle \xi \rangle^{\mu/2},
\end{equation}
where $ N $ is a large integer and $\mu$ is any real
number. Then we may define the class $ S(m_{N}^{t, \mu}, g) $ of symbols in the standard way.
We point out explicitly that $t$ is just a parameter and at this level we may omit it in our notation.
We have

\begin{prop}
$ f^{-N}(t, \xi) \in S(m_{N}^{t, 0}, g) $.
\end{prop}

We have also

\begin{prop}
Let $ c(x, \xi) \in S^{\mu}(1, g) $ be a classical symbol of order $\mu$. Then 
$ f^{-N}(t, \xi) \#_{x} c (x, \xi)  =  b_{t}(x, \xi),$
where $ b_{t} \in S(m_{N}^{t,\mu}, g)$. Here $ \#_{x} $ denotes the
operation of formal asymptotic composition $ g \#_{x} c =\sum_{|\alpha| \geq 0} \frac{1}{\alpha!}\partial_{\xi}^{\alpha} g (x, \xi)D_{x}^{\alpha}c (x, \xi) $.
\end{prop}

To examine the lower order terms we need to handle the term

$$f^{-N}(t, D_x) b(\epsilon^{2/3}t, \epsilon x, D_x)f^N(t, D_x),$$
$b(t, x, D_x)$ being a second order pseudodifferential operator. We deduce that 
$$B_N = f^{-N}(t, \xi)\#_x b(\epsilon^{2/3} t, \epsilon x, \xi)(1 + |\xi|^2)^{-1} f^N(t, \xi)\in S(m_N^{t, 0}, g)$$
but we need to estimate the $L^2$ norm of the operator $B_N$ and for this reason we take $\epsilon$ to be of order
${\mathcal O}(\frac{1}{N})$.  Therefore in the calculus of lower order terms of $B_N$ the powers of $N$ are compensated by the powers of $\epsilon$. Moreover, we may write the composition of symbols $B_N$ by using a finite sum and an integral representation of the remainder introduced by J.M. Bony \cite{Bo}.\\

 The details of the analysis of the operators with variable coefficients depending on $(t, x)$ will be given in a paper in preparation \cite{BBP}.

\end{document}